\newcommand{\bull}{\vrule height .9ex width .8ex depth -.1ex}
 \newcommand{\ppp}{\hfill $\bull$ }

 \documentclass[11pt]{article}
 \author{ Mohammed-Larbi Labbi}
   \title{On  Weitzenb\"ock Curvature Operators }
   \date{}
   \usepackage{amsmath}
   \usepackage{amscd}
\newtheorem{theorem}{Theorem}[section]
\newtheorem{corollary}[theorem]{Corollary}

\newtheorem{proposition}[theorem]{Proposition}

\begin{document}
   \maketitle
   \begin{abstract} 
The Weitzenb\"ock curvature operators are the curvature terms of order zero 
 that appear in the well known
classical Weitzenb\"ock formula. In this paper, we use the formalism of double forms
to prove a simple formula for this operators and to study their geometric properties.
   \end{abstract}
 \section{Introduction}
The Weitzenb\"ock curvature operator, denoted throughout this paper by ${\cal N}$,
is the curvature term of order zero (i.e. depends
linearly on the Riemann curvature tensor)
 that appears in the well known classical
Weitzenb\"ock formula. The former expresses the Laplacian
$\Delta$ of differential forms in terms of the
  Levi-Civita connexion $\nabla$, precisely we have:\\
$$\Delta=\nabla^*\nabla+{\cal N}.$$
This formula is particularly important in the study of interactions between the geometry and 
topology of  manifolds. In fact, there exists a method, 
due to Bochner and known as vanishing theorems,
consisting of proving the vanishing of Betti numbers of a Riemannian manifold with
 a positive curvature condition stronger than the positivity of  
 ${\cal N}$. This method mainly applies to compact manifolds.\\
For each $p$, the curvature operator 
 ${\cal N}$ preserves $p$-forms and it is self-adjoint. Therefore, by duality, we can consider
it as a double form.\\
The main obstacle with  ${\cal N}$ is that it is actually a complicated expression of 
the curvature. Several  simplifications of such expression exist in the literature.  For instance, 
the Clifford
formalism of Lawson-Michelsohn \cite{Lawson-Michelsohn}, see also \cite{Labbi-betti},
the works of  Gallot-Meyer \cite{Gallot-Meyer}, Maillot \cite{Maillot}, Bourguignon
\cite{Bourguignon1}  ... \\
In this paper, using the formalism of double forms, we prove a simple formula for
${\cal N}$ and then we use it to study some geometric properties of this curvature.
This formula was first established by Bourguignon, see proposition 8.6 in \cite{Bourguignon1}.
Unfortunately, up to my knowledge, this nice formula is not well known  
and it is not used even though the paper  \cite{Bourguignon1} is very famous!.
 From my side, I noticed the existence of this formula  only once
my proof was finished. However our proof is completely different. 
It is of algebraic nature and direct.\\

\section{Double Forms}
Let $(V,g)$ be an Euclidean real vector space  of dimension n. In the
 following
we shall identify whenever convenient (via their Euclidean structures),
the vector spaces
 with their duals. Let
  $\Lambda V=\bigoplus_{p\geq 0}\Lambda^{p}V$ denotes the exterior algebra
 of  $p$-vectors on $V$. \\
A double form on $V$ of degree $(p,q)$ can be defined as a 
 bilinear form $\Lambda^pV\times\Lambda^qV\rightarrow {\bf R}$. That is
     a multilinear form which is skew symmetric in the first $p$-arguments and also
     in the last $q$-arguments. If $p=q$ and the bilinear form is symmetric we say that we
 have a symmetric double form.\\
Double forms are abundant in geometry. The Riemann curvature tensor and the Weyl curvature tensor are
 symmetric double forms
of order $(2,2)$. The metric, Ricci, Einstein, and Schouten tensors are symmetric double forms
of degree $(1,1)$. Gauss-Kronecker tensors \cite{Thorpe} and the Weitzenb\"ock curvatures 
(see the next section)
are examples of symmetric double forms of higher order.\\
The exterior product between $p$-vectors extends in a natural way to double forms of any degree
 and we obtain
the so called
Kulkarni-Nomizu product. For a $(p,q)$-double form $\omega_1$ and an $(r,s)$-double form
$\omega_2$, the product $\omega_1\omega_2$ is a double form of degree $(p+r,q+s)$ given by
 \begin{equation}
\begin{split}
\label{eps:prod}
 &\omega_1.\omega_2(x_1\wedge...\wedge x_{p+r},y_1\wedge...\wedge y_{q+s})\\
&= {1\over p!r!s!q!}\sum_{\sigma\in S_{p+r}, \rho\in S_{q+s}}
\epsilon(\sigma)\epsilon(\rho)
\omega_1(x_{\sigma(1)}\wedge...\wedge x_{\sigma(p)},y_{\rho(1)}
\wedge...\wedge y_{\rho(q)})\\
&\phantom{...mmmmmmmmmmm}
\omega_2(x_{\sigma(p+1)}\wedge...\wedge x_{\sigma(p+r)},y_{\rho(q+1)}
\wedge...\wedge y_{\rho(q+s)})
\end{split}
\end{equation}  
In particular , the product of the inner product  $g$ with itself
 $k$-times determines the canonical scalar product on
 $\Lambda^{p}V$. Precisely we have
$$g^k(x_1 \wedge...\wedge x_k,y_1\wedge...\wedge y_k)=k!\det[g(x_i,y_j)].$$
Let us denote by $S^2\Lambda^pV$ the set of all the symmetric double forms of order $(p,p)$
on $V$ and let $S^2\Lambda V = \bigoplus_{p=0}^{p=n}S^2\Lambda^pV$.
 Recall that it is a commutative algebra.
\subsection{The multiplication map by $g^k$ and the contraction map}
The multiplication map by the powers of the metric, that is $g^k$, in $S^2\Lambda^pV$
plays an important role in the study of double forms.
In \cite{Labbi-double} we proved the following fundamental  property of this map:
\begin{proposition}[\cite{Labbi-double}]\label{mult:injec} 
Let $\omega_1,\omega_2$ be two $(p,p)$-double forms  and $k\leq n-2p$, then
\begin{equation}
g^k\omega_1=g^k\omega_2\Longrightarrow \omega_1=\omega_2.
\end{equation}
\end{proposition}
Note here for future reference that, for each $p$, the natural scalar product $\frac{g^p}{p!}$ on   $\Lambda^{p}V$ induces canonically
 a natural scalar product on $S^2\Lambda^pV$, we shall denote it by $ \langle,\rangle$.
 We extend  $\langle,\rangle$ to  $\bigoplus_{p=0}^{p=n}S^2\Lambda^pV$ by declaring that   $S^2\Lambda^pV \perp 
S^2\Lambda^qV$
if $p\not = q$.
\\
A second fundamental map on double forms is the contraction map $c$. It decreases the order of a
double form by $1$ and it is
 the adjoint of the multiplication map by $g$, precisely we have.
\begin{proposition}[\cite{Labbi-double}]\label{theo:gc} For arbitrary double forms
 $\omega_1, \omega_2$ we have
\begin{equation}
\label{adj:gc}
\langle g\omega_1,\omega_2\rangle=\langle\omega_1,c\omega_2\rangle.
\end{equation}
In particular, for every $k\geq 1$, we have
$\langle g^k\omega_1,\omega_2\rangle=\langle\omega_1,c^k\omega_2\rangle$. Where $c^k$ denotes the composition
of the contraction map $c$ with itself $k$-times.
\end{proposition}

\subsection{The first Bianchi identity and sectional curvatures}
We say that a double form $\omega$ of degree $(p,q)$, with $q\geq 1$,
 satisfies the first Bianchi identity if for all vectors $(x_i),(y_j)$ we have
$$\sum_{j=1}^{p+1}(-1)^j\omega(x_1\wedge...\wedge \hat{x}_j\wedge ... x_{p+1},
x_j\wedge y_1\wedge...\wedge y_{q-1})=0,$$
where $\hat{}$ denotes omission. Note that this identity is preserved by multiplication of double
forms.
Let us then denote by $S^2_1(\Lambda V)$ the sub-algebra of $S^2(\Lambda V)$ consisting of all
symmetric double forms  satisfying the first Bianchi identity. This space is not
irreducible under the natural action of the orthogonal group. In fact, Kulkarni \cite{Kulkarni}
proved that the full reduction into irreducible components is given by
\begin{equation}\label{decomposition}
S^2_1(\Lambda^pV)=E^p_1\oplus gE^{p-1}_1\oplus  g^2E^{p-2}_1\oplus ...\oplus g^pE^{0}_1.
\end{equation}
Where, for each $0\leq k\leq p$, $E_1^k=\{\omega\in S_1^2(\Lambda^kV) :c\omega =0\}$.\\
Another important property of double forms satisfying the first Bianchi identity is that they
are determined by their sectional curvatures. Recall that the sectional curvature of a given
symmetric $(p,p)$-double form $\omega$ is a function, say $K_\omega$, defined on the
Grassman algebra of $p$-planes in $V$. For a $p$-plane $P$, we set 
 $$K_\omega(P)=\omega(e_1\wedge...\wedge e_p,e_1\wedge...\wedge e_p),$$
 where $\{e_1,...,e_p\}$ is any orthonormal basis of $P$. 
We have the following characterization of symmetric double forms satisfying the first Bianchi identity
and with constant sectional curvature \cite{Kulkarni}:
  \begin{equation}\label{kom:const}
  K_{\omega}\equiv c \quad \text{is constant}\quad \text{ if and only if}\qquad
  \omega=c{g^p\over p!}.
  \end{equation}
 For  a symmetric $(r,r)$-double form $\omega$ and every $p$, $0\leq p\leq n-r$,
we can use  formula (\ref{eps:prod}) to evaluate the sectional curvature of the products
 $g^p\omega$, this will be used later in this paper.
 Let 
 $\{e_1,...,e_{p+r}\}$ be orthonormal vectors, then
 \begin{equation}\label{gpom:form}
 \begin{split}
  g^p\omega&(e_1\wedge...\wedge e_{p+r},e_1\wedge...\wedge e_{p+r})\\
 =&p!\sum_{1\leq i_1<i_2<...<i_r\leq p+r}\omega(e_{i_1}\wedge...\wedge e_{i_r},
 e_{i_1}\wedge...\wedge e_{i_r}).
 \end{split}
 \end{equation}

\subsection{Generalized Hodge star operator}
We suppose here that an orientation is fixed on the vector space  $V$.
The classical
Hodge star operator  $*:\Lambda^{p}V\rightarrow \Lambda^{n-p}V$ can be extended in a 
natural way to double forms. For a double form $\omega$, set 
 $*\omega(.,.)=\omega(*.,*.)$. Many classical  properties can be generalized 
to the case of this new operator, see \cite{Labbi-double}. 
The generalized Hodge star operator was shown to be an important tool in the study of 
double forms. In particular, it provides a second simple relation between the contraction map
$c$ and the multiplication map by  $g$, as follows:\\
\begin{proposition}[\cite{Labbi-double}] For every $(p,p)$-double form
 $\omega$, we have
\begin{equation}
\label{gstar:c}
g\omega=*c*\omega.
\end{equation}
 In particular, for every  $k\geq 1$, we have $g^k\omega=*c^k*\omega$.
\end{proposition}

\section{Weitzenb\"ock Operators}
\subsection{Clifford multiplication of $p$-vectors}
Let $e\in V$,
recall that the interior product $i_e$ on $\Lambda V$ is the adjoint  of 
the exterior multiplication map by $e$.
We define the Clifford multiplication of two $p$-vectors, denoted by a dot, as follows:
For $e\in V$ and $\omega\in \Lambda V$, set
$$e.\omega=e\wedge \omega-i_e\omega.$$
In particular, for $e,f\in V$ we have $e.f=e\wedge f-g(e,f)$.\\
Assuming this product to be associative, we can define  the products 
$ e_{i_1}..... e_{i_{p}}$ for all $p$ and $ e_{i_k}\in V$. Then  using linearity, we can extend it
 to
a product on $\Lambda V$: the Clifford product. \\
It is not difficult to check that we can recover the exterior product as:
\begin{equation}
 e_{i_{1}}\wedge...\wedge e_{i_{p}}=\frac{1}{p!}
\sum_{\sigma\in {S_p}}\epsilon(\sigma) e_{\sigma(i_{1})}.....\wedge e_{\sigma(i_{p})}.
\end{equation}
\subsection{Definition of Weitzenb\"ock operators}
For $\phi \in \Lambda^2V$, we define the linear operator $ad_\phi$ on $\Lambda V$ by
$$ad_\phi(\psi)=[\phi,\psi]=\phi.\psi-\psi.\phi.$$
In particular, a straightforward computation shows that
 for orthonormal vectors $\{e_i\}$ we have
\begin{equation}\label{ad}
ad_{e_i.e_j}e_{i_{1}}..... e_{i_{p}}=\left\{
\begin{array}{ccccc}
0 &\text{if}& i,j\in\{i_1,...,i_p\} &\text{or}&i,j\not\in \{i_1,...,i_p\}\\
2e_i.e_j.e_{i_{1}}..... e_{i_{p}}, &\, &\text{otherwise.}&\, &\\
\end{array}
\right.
\end{equation}
Let now $\omega$ be a symmetric $(2,2)$-double form. We define the Weitzenb\"ock transformation
of order $p$ at $\omega$ to be the symmetric $(p,p)$-double form ${\cal N}_p(\omega)$
defined as follows, see \cite{Lawson-Michelsohn, Labbi-these, Labbi-betti}.
\begin{equation}\label{Ndefinition}
{\cal N}_p(\omega)(\psi_1,\psi_2)=\frac{1}{4}\sum_{\scriptstyle i<j, k<l}
\omega(e_i\wedge e_j,e_k\wedge e_l)\langle ad_{e_i.e_j}\psi_1,ad_{e_k.e_l}\psi_2\rangle.
\end{equation}
Where $(e_1,...,e_n)$ denotes an arbitrary orthonormal basis of $V$.\\
This definition is of course motivated by the curvature term in the Weitzenb\"ock formula.
 Note that ${\cal N}_p$ is a linear operator.\\
\subsection{Sectional curvatures of Weitzenb\"ock}
Let $P$ be a $p$-plane in $V$ spanned by orthonormal vectors $\{e_1,...,e_p\}$.
A direct computation using formula (\ref{ad}) shows that
\begin{equation}\label{sectional}
\begin{split}
{\cal N}_p&(\omega)(e_1\wedge ...\wedge e_p,e_1\wedge ...\wedge e_p)\\
&=\frac{1}{4} \sum_{ k<l}\sum_{i=1}^p\sum_{j=p+1}^n
\omega(e_i\wedge e_j,e_k\wedge e_l)(-1)^{i+1}\langle 2e_j.e_1...\hat{e_i}...e_p,ad_{e_k.e_l}e_1...e_p\rangle\\
&=\frac{1}{4} \sum_{i=1}^p\sum_{j=p+1}^n
\omega(e_i\wedge e_j,e_i\wedge e_j)\langle 2e_j.e_1...\hat{e_i}...e_p,2e_j.e_1...\hat{e_i}...e_p\rangle\\
&=\sum_{i=1}^p\sum_{j=p+1}^n\omega(e_i\wedge e_j,e_i\wedge e_j).\\
\end{split}
\end{equation}
Now  using (\ref{gpom:form}) we get for $p\geq 2$:
\begin{equation}
\begin{split}
\sum_{i=1}^p\sum_{j=p+1}^n\omega(e_i\wedge e_j &,e_i\wedge e_j)=\sum_{i=1}^p\left\{
c\omega(e_i,e_i)
-\sum_{i,j=1}^p\omega(e_i\wedge e_j,e_i\wedge e_j)\right\}\\
&=\{\frac{g^{p-1}}{(p-1)!}c\omega -2\frac{g^{p-2}}{(p-2)!}
\omega\}(e_{1}\wedge...\wedge e_{p},
e_{1}\wedge...\wedge e_{p}).
\end{split}
\end{equation}
We conclude  that the two double forms
 ${\cal N}_p$ and 
$\{\frac{g^{p-1}}{(p-1)!}c\omega -2\frac{g^{p-2}}{(p-2)!}\omega\}$ have the same sectional
curvatures for every $2\leq p\leq n-2$. It is then natural to expect the equality of these two
double forms.
This is in fact true and shall be proved  in the section below.

\subsection{A simple formula for the Weitzenb\"ock Transformation}
In the proposition below, we prove a simple formula for the transformation ${\cal N}$.
\begin{theorem}
Let $\omega$ be a symmetric $(2,2)$-double form satisfying the first Bianchi identity, then
for every $2\leq p\leq n-2$, the Weitzenb\"ock transformation of order $p$ is determined by
\begin{equation}\label{main-equation}
{\cal N}_p(\omega)=\left\{\frac{gc\omega}{p-1}-2\omega\right\}\frac{g^{p-2}}{(p-2)!}.
\end{equation}
\end{theorem}
{\sc Proof.} We shall prove the equality by evaluating both terms of the equation on decomposed
$p$-vectors. 
Let
$\{e_1,...,e_{n}\}$ be an arbitrary orthonormal basis of $V$ and
$p\geq 1$. Let
$e_I=e_{i_1}\wedge...\wedge e_{i_{p}}$ and
$e_J=e_{j_1}\wedge...\wedge e_{j_{p}}$ be two arbitrary
elements of the standard basis of $\Lambda^{p} V$.\\
First,
it is straightforward from formulas (\ref{eps:prod}) and  (\ref{Ndefinition}) that 
if $\text{card}(I\cap J)\leq p-3$ then 
$${\cal N}_p(\omega)(e_I,e_J)=\left\{\frac{gc\omega}{p-1}-2\omega\right\}\frac{g^{p-2}}{(p-2)!}(e_I,e_J)=0.
$$\\
Next, suppose $\text{card}(I\cap J)=p-2$. Without loss of generality,
we may re-index as follows: 
$e_I=e_1\wedge...\wedge e_p$ and $e_J=e_1\wedge...\wedge e_{p-2}
\wedge e_{p+1}\wedge e_{p+2}$.
  Let  
$f_I=f_1\wedge...\wedge f_p=e_J$, where
\[ \left\{
\begin{array}{ccc}
f_j=e_j& \mbox{if}& j<p-1\\
f_{p-1}=e_{p+1}& \mbox{and} & f_p=e_{p+2}\\
\end{array}
\right. \]
Let $\omega_1=\frac{1}{(p-2)!} \{\frac{gc\omega}{(p-1)}-2\omega\}$,
 then using  formula   (\ref{eps:prod}) for the product, the right hand side 
of  equation
(\ref{main-equation})
is given by:
\begin{equation}\label{produit}
\begin{split}
g^{p-2}\omega_1(e_I,e_J)=
 \sum_{\sigma ,\rho \in Sh(2,p-2)}
\epsilon(\sigma)&\epsilon(\rho)
\omega_1(e_{\sigma(1)}\wedge e_{\sigma(2)},f_{\rho(1)}
\wedge f_{\rho(2)})\\
&g^{p-2}(e_{\sigma(3)}\wedge...\wedge e_{\sigma(p)},f_{\rho(3)}
\wedge...\wedge f_{\rho(p)}).
\end{split}
\end{equation}
Where $\sigma\in Sh(2,p)$ means a   permutation of $\{1,...,p\}$
such that $\sigma(1)< \sigma(2)$
and $\sigma(3)<...<\sigma(p)$.\\
Remark that the only permutations which yield non zero terms in
the previous summation are
$$ \sigma=\rho=\left(
\begin{array}{ccccc}
1&2&3&...& p\\
p-1&p&1&...&p-2 \\
\end{array}
\right) $$
Since $\epsilon(\sigma)=\epsilon(\rho)=1$  we get
\begin{equation*}
\begin{split}
g^{p-2}\omega_1(e_I,e_J)&=(p-2)!\omega_1(e_{p-1}\wedge e_p,e_{p+1}\wedge e_{p+2})\\
&=-2\omega(e_{p-1}\wedge e_p,e_{p+1}\wedge e_{p+2}).\\
\end{split}
\end{equation*}
On the other hand, using (\ref{ad}) we have
\begin{equation}\label{wsum}\begin{split}
4{\cal N}_p(\omega)(e_I,e_J)=\sum_{k<l}\sum_{i=1}^p\sum_{j=p+1}^n
&\omega(e_i\wedge e_j,e_k\wedge e_l)(-1)^{i+1}
\\
&\langle 2e_j.e_1...\hat{e_i}....e_p,
 ad_{{e_k}.{e_l}}e_1.....e_{p-2}e_{p+1}e_{p+2}\rangle.
\end{split}
\end{equation}
Note that the terms with $i\leq p-2$ vanish. For $i=p-1$ (resp. for $i=p$) , the only non zero terms are
those corresponding to $k=p,l=p+1,j=p+2$ or $k=p,l=p+2,j=p+1$ (resp. 
$k=p-1,l=p+1,j=p+2$ or $k=p-1,l=p+2,j=p+1$). Therefore the previous summation reduces to
\begin{equation}
\begin{split}
4{\cal N}_p(\omega)(e_I,e_J)&=8\omega(e_{p-1}\wedge e_{p+2},e_{p}\wedge e_{p+1})
-8\omega(e_{p-1}\wedge e_{p+1},e_{p}\wedge e_{p+2})\\
&=8\omega(e_{p-1}\wedge e_{p},e_{p+2}\wedge e_{p+1}).
\end{split}
\end{equation}
Where in the last equality we have used the the first Bianchi identity.\\
We now study the next case  where
$\text{card}(I\cap J)=p-1$. Here also without loss of generality
we may assume  
$e_I=e_1\wedge...\wedge e_p$ and $e_J=e_1\wedge...\wedge e_{p-1}
\wedge e_{p+1}$.
  Let  
$f_I=f_1\wedge...\wedge f_p=e_J$, where
\[ \left\{
\begin{array}{ccc}
f_j=e_j& \text{if}& j<p\\
f_p=e_{p+1}.&\, &\\
\end{array}
\right. \]
The product $g^{p-2}\omega_1(e_I,e_J)$ has the form (\ref{produit}).
 In this case,
the only permutations which
 yield non zero terms are the following:\\
\[ \sigma_k=\rho_k=\left(
\begin{array}{ccccccc}
1&2&3&...&k&...&p\\
k&p&1&...&\hat{k} &...&p-1\\
\end{array}
\right),
\]
where $k$ ranges from $1$ to $p-1$.
Therefore we immediately obtain
\begin{equation}
\begin{split}
g^{p-2}\omega_1(e_I,e_J)&=(p-2)!\sum_{k=1}^{p-1}\omega_1(e_k\wedge e_p,e_k \wedge e_{p+1})\\
&=\sum_{k=1}^{p-1}\left\{ \frac{1}{p-1}c\omega(e_p,e_{p+1})-2\omega(e_k\wedge e_p,e_k 
\wedge e_{p+1})\right\}\\
&=c\omega(e_p,e_{p+1})-2\sum_{k=1}^{p-1}\omega(e_k\wedge e_p,e_k \wedge e_{p+1}).\\
\end{split}
\end{equation}
We now evaluate the right hand  term on the same $p$-vectors.
Recall that ${\cal N}_p(\omega)(e_I,e_J)$ can be written as a summation as in (\ref{wsum}). 
In this summation and  for $i\leq p-1$  (resp. for $i=p$), the only non zero terms are
those corresponding to $k=i,l=p,j=p+1$  (resp. 
$k=p+1,l=j$). Therefore we can write
\begin{equation}
\begin{split}
4{\cal N}_p(\omega)(e_I,e_J)&=-4\sum_{i=1}^{p-1}
\omega(e_{i}\wedge e_{p+1},e_{i}\wedge e_{p})
+4\sum_{j=p+1}^n\omega(e_{i}\wedge e_{p+1},e_{i}\wedge e_{p})\\
&=4c\omega(e_p,e_{p+1})-8\sum_{i=1}^{p-1}\omega(e_i\wedge e_{p+1},e_i \wedge e_{p}).\\
\end{split}
\end{equation}
Finally, In the case where $\text{card}(I\cap J)=p$, then $I=J$. This case was  already 
discussed in the previous section (that is the two double forms
 have the same sectional curvatures). 
This completes the proof of the theorem. \ppp
\par\medskip\noindent
As a consequence of the previous theorem, the Weitzenb\"ock transformation preserves the
 first Bianchi identity. That  is, it sends  symmetric
$(2,2)$-double forms satisfying the first Bianchi identity to symmetric $(p,p)$-double
forms satisfying the first Bianchi identity. Furthermore we have:
\begin{corollary}
For every $2\leq p\leq n$ the Weitzenb\"ock linear transformation
$${\cal N}_p:S_1^2(\Lambda^2V)\rightarrow S_1^2(\Lambda^pV)$$
 is injective. Furthermore, its adjoint operator  is given by
$${\cal N}_p^*\omega=\left(\frac{gc^{p-1}}{(p-1)!}-2\frac{c^{p-2}}{(p-2)!}
\right)\omega$$
\end{corollary}
{\sc Proof.} The first claim results directly from the previous theorem and 
 proposition \ref{mult:injec}. The second claim  is a direct consequence of
proposition \ref{theo:gc}.\ppp

\section{Geometric Properties}
In this section, Let  $\omega$ be a fixed symmetric $(2,2)$-double form satisfying the first
Bianchi identity on an $n$-dimensional
vector space $V$ (keeping in mind the typical example of the Riemann
curvature tensor of a Riemannian $n$-manifold at a given point).
In order to simplify the notations we shall denote by ${\cal N}_p$ 
the double form
${\cal N}_p(\omega)$.\\
The following result can be proved easily using the previous theorem after noticing that both double forms satisfy
the first Bianchi identity and they have the same sectional curvatures (see formula
\ref{sectional}) :
\begin{proposition}
For each $p$ such that $ 2\leq p\leq n-2$, we have
\begin{equation}\label{np}
* {\cal N}_p={\cal N}_{n-p}.
\end{equation}
In particular, if $n=2p$ then $*{\cal N}_p={\cal N}_{p}.$
\end{proposition}
It is clear from the previous theorem  that the double form  ${\cal N}_p$ is divisible by $g^{p-2}$ and hence its
orthogonal decomposition is reduced as follows:

\begin{proposition}
With respect to the irreducible orthogonal decomposition  (\ref{decomposition}),
if $\omega=\omega_2+g\omega_1+g^2\omega_0$, then 
the double form  ${\cal N}_p$, $2\leq p\leq n-2$,  splits as follows: 
$${\cal N}_p=g^{p-2}\left\{ \frac{-2\omega_2}{(p-2)!}\right\}
+g^{p-1}\left\{\frac{(n-2p)\omega_1}{(p-1)!}\right\}+
g^p\left\{\frac{2(n-p)\omega_0}{(p-1)!}\right\}.$$
\end{proposition}
Let now $\omega=R$ be the Riemann curvature tensor of a given Riemannian manifold. 
Recall that the sectional curvature of Weitzenb\"ock is given by (\ref{sectional})
with $\omega=R$, that is a generalization of the Ricci curvature.
Also, a double form of order  $(p,p)$ and satisfying the first Bianchi identity is with
constant sectional curvature if and only if it is proportional to the metric $g^p$.
The following corollary is therefore  a straightforward consequence of the previous proposition. 
\begin{corollary} For $2\leq p\leq n-2$, a Riemannian manifold  $(M,g)$
of  dimension
 $n$ has its  Weitzenb\"ock  
sectional curvature of order  $p$ constant if and only if it is either with constant sectional
curvature or conformally flat with dimension  $n=2p$.
\end{corollary}
Let us note here that the previous corollary was first proved by Tachibana in  \cite{Tachibana}. Where  the sectional curvature of Weitzenb\"ock  is called the mean curvature of a $p$-plane.\\
We now go back to the general situation where $\omega$ is any algebraic $(p,p)$-double form.
Using formula (12) in \cite{Labbi-double} and the previous proposition, a long but 
straightforward
computation shows that:
\begin{theorem}
For every $2\leq p\leq n-2$, the contractions of the double form
${\cal N}_p$ up to the order
$p$ and $p-1$ are respectively given by
\begin{equation}\label{contractions}
\begin{split}
c^p({\cal N}_p)&=\frac{p(n-2)!}{(n-p-1)!}c^2\omega,\\
c^{p-1}({\cal N}_p)&=\frac{(n-3)!}{(n-p-1)!}\left\{ (n-2p)c\omega+(p-1)c^2\omega g\right\}.\\
\end{split}
\end{equation}
Furthermore, for any $k\leq p-2$, the contraction up to the  order $k$
 is given by
\begin{equation}\label{contractions2}
\begin{split}
c^{k}({\cal N}_p)=\frac{(n-p+k-2)!g^{p-k-2}}{(n-p-2)!(p-k-2)!}\{
-2\omega &+\frac{n-k-p-1}{(n-p-1)(p-k-1)}gc\omega\\
&+
\frac{k}{(n-p-1)(p-k-1)(p-k)}g^2c^2\omega\}.
\end{split}
\end{equation}
\end{theorem}
{\sc Remark.} Let $E=\frac{1}{2}c^2\omega g-c\omega$ be the Einstein curvature tensor of $\omega$, then the
second  formula in 
(\ref{contractions}) may be written in the following alternative form:
\begin{equation}\label{alternative}
c^{p-1}({\cal N}_p)=\frac{(n-3)!}{(n-p-1)!}\left\{ \frac{n-2}{2}c^2\omega g
-(n-2p)E \right\}.
\end{equation}
{\sc Proof.} Let $k\leq p-2$, then using formulas (\ref{adj:gc}),(\ref{np})  we get
\begin{equation*}
\begin{split}
c^{k}({\cal N}_p)=*g^k*({\cal N}_p)&=*g^k{\cal N}_{n-p}\\
&=*\left\{(\frac{1}{n-p-1}c\omega.g-2\omega)\frac{g^{n-p-2+k}}{(n-p-2)!}\right\}.
\end{split}
\end{equation*}
A direct application of  formula (15) in \cite{Labbi-double} shows then that
\begin{equation*}
\begin{split}
c^{k}({\cal N}_p)&=-\frac{(n-p-1+k)!}{(n-p-1)!}\frac{g^{p-k-1}}{(p-k-1)!}c\omega
+\frac{(n-p-1+k)!}{(n-p-1)!}\frac{g^{p-k}}{(p-k)!}c^2\omega\\
&-2\frac{(n-p-2+k)!}{(n-p-2)!}\frac{g^{p-k-2}}{(p-k-2)!}\omega
+2\frac{(n-p-2+k)!}{(n-p-2)!}\frac{g^{p-k-1}}{(p-k-1)!}c\omega\\
&-\frac{(n-p-2+k)!}{(n-p-2)!}\frac{g^{p-k}}{(p-k)!}c^2\omega.
\end{split}
\end{equation*}
Note that for $k=p$ (resp. $k=p-1$) the first, third and fourth terms should be dropped
from the previous summation
 (resp. the third term). A direct computation then
 yields the desired results.\ppp 

\section{Positivity of the Weitzenb\"ock curvatures}
Let now $\omega=R$ be the curvature tensor of a Riemannian manifold and  ${\cal N}_p$ be 
the corresponding
Weitzenb\"ock double form.\\
The positivity of  ${\cal N}_p$ is clearly very important in the study of the interactions 
between the geometry and topology of manifolds.
Therefore it is particularly interesting to understand this condition and to relate it
to the positivity of the other well known curvatures.\\
The case $p=1$ (and by duality $p=n-1$) is well known since ${\cal N}_1$ coincides
 with the Ricci curvature, so we suppose $2\leq p\leq n-2$.\\
Needless to say at this stage  that the positivity of
 ${\cal N}_p$ is strictly stronger than the positivity of its sectional curvature
 (likewise the positivity of the Riemann curvature operator and of its sectional curvature).
 We already know that the positivity of the Riemann curvature operator implies the one of  
${\cal N}_p$ for every $p$. This is a theorem of  Meyer \cite{Gallot-Meyer}, see 
 \cite{Lawson-Michelsohn} for a simplified proof of this result. \\
On the other side, it is clear from equation (\ref{contractions}) that the positivity
of ${\cal N}_p$ implies the positivity of the scalar curvature. Consequently, the positivity
of the Weitzenb\"ock curvatures are intermediate between positive scalar curvature and
positive curvature operator.\par\medskip\noindent
The following proposition is a direct consequence of formulas (\ref{contractions}),
(\ref{contractions2}) and (\ref{alternative}):
\begin{proposition} For each $p$ with $0\leq p\leq n-1$,
a  Riemannian $n$-manifold  with positive scalar curvature has its contracted
Weitzenb\"ock curvature up to the order $p$, that is $c^{p}({\cal N}_p)$, positive.\\
Furthermore, in each of the following cases the contracted curvature $c^{p-1}({\cal N}_p)$ is positive:
\begin{enumerate}
\item If $n=2p+2$ and  $(M,g)$ has positive scalar curvature.
\item If $n\leq 2p+2$ and  $(M,g)$ has positive Einstein
 curvature.
\item If $n\geq 2p+2$ and  $(M,g)$ has positive Ricci
 curvature.
\end{enumerate}  
 \end{proposition}
In what follows, we investigate the relation to the positivity of the $p$-curvature 
\cite{Labbi-betti, Labbi-these}. Recall that the $p$-curvature $s_p$ of the Riemann 
curvature tensor $R$ is defined to be  the sectional curvature 
of the tensor $*\left( \frac{1}{(n-p-2)!}g^{n-p-2}R \right)$. \\
It is clear that the $p$-curvature appears naturally in the formula (\ref{main-equation})
for ${\cal N}_p$
but it contributes with a negative sign!. The influence of the positivity of the
$p$-curvature on the positivity of ${\cal N}_p$ becomes clear in the following formula:
\begin{proposition} For every $2\leq p\leq n-2$ such that $n+p$ is even we have
\begin{equation}
{\cal N}_{\frac{n+p}{2}}=C(n,p) g^{\frac{n-p}{2}}
\left\{*(\frac{p(p-1)}{(n-p-2)!}g^{n-p-2}R)-\frac{(n-1)(n-2)}{(p-2)!}
g^{p-2}W\right\}.
\end{equation}
Where $C(n,p)=\frac{2(p-2)!}{(\frac{n+p-4}{2})!(n+p-2)(n-p-1)}$ and $W$ is the double form 
associated to the standard weyl curvature tensor.
\end{proposition}
{\sc Proof.} Let $R=\omega_2 +g\omega_1+g^2\omega_0$ be the standard decomposition of the 
Riemann curvature tensor, where $\omega_2=W$. Applying formula (20) in \cite{Labbi-double} we obtain
$$*\frac{g^{n-p-2}R}{(n-p-2)!}=\frac{g^{p-2}}{(p-2)!}\left\{\omega_2-\frac{n-p-1}{p-1}g\omega_1
+\frac{(n-p)(n-p-1)}{p(p-1)}g^2\omega_0\right\}.$$
Inserting this in the right hand side of the desired equation one can  easily  
reach to the left hand side. \ppp \\
As a direct consequence of the previous formula we get an alternative proof of our result
in \cite{Labbi-betti}  on the vanishing of betti numbers around the middle dimension
for conformally flat manifolds ($W=0$) with positive $p$-curvature.\\
Finally, taking into consideration the existing similarity between positive isotropic curvature and positive $p$-curvature, see 
\cite{Labbi-isotropic, Labbi-surgery}, the previous proposition may help in answering the following question:\\
\begin{center}
For a compact $n$-dimensional Riemannian manifold, does positive isotropic curvature imply the vanishing of the betti numbers $b_i$ for
$2\leq i\leq n-2$?
\end{center}
With reference to the previous proposition, a positive answer to the  question above  would be possible if one can prove that the  positivity of isotropic curvature implies the positivity of  the operator
$$*(\frac{p(p-1)}{(n-p-2)!}g^{n-p-2}R)-\frac{(n-1)(n-2)}{(p-2)!}
g^{p-2}W,$$
for all $0\leq p\leq n-4$.

\vspace{1cm}
\noindent
Labbi Mohammed-Larbi\\
  Department of Mathematics,\\
 College of Science, University of Bahrain,\\
  P. O. Box 32038 Bahrain.\\
  E-mail: labbi@sci.uob.bh


\begin{thebibliography}{99}


\bibitem{Bourguignon1}   Bourguignon, J. P. \emph{Les vari\'et\'es de dimension $4$ \`a 
signature non nulle dont la courbure est harmonique sont d'Einstein}, Invent. Math. 63, 263-286
(1981).


\bibitem{Gallot-Meyer} Gallot, S., Meyer, D.  \emph{Op\'erateur de courbure et Laplacien des 
formes diff\'er\'erentielles d'une vari\'et\'e riemannienne,} J. Math. pures et appl., 54,
259-284 (1975).


\bibitem{Kulkarni} Kulkarni, R. S.,
\emph{On Bianchi Identities}, Math. Ann. 199, 175-204 (1972).





\bibitem{Labbi-betti} Labbi, M.-L., \emph{ Sur les nombres de Betti
des vari\'et\'es conform\'ements
plates}, CRAS, t 319, s\'erie I, 77-80 (1994).

\bibitem{Labbi-these} Labbi, M.-L., \emph{ Vari\'et\'es riemanniennes \`a $p$-courbure
positive,} th\`ese, publication of Montpellier II University (1995), France.

\bibitem{Labbi-isotropic} Labbi, M.-L., \emph{ On compact manifolds with positive 
isotropic curvature},
Proceedings of the American Mathematical Society, Volume 128, Number 5, Pages 1467-1474 (2000).

\bibitem{Labbi-surgery} Labbi, M.-L., \emph{On positive 
isotropic curvature and surgeries}, Journal
of Differential Geometry and its applications, 17, 37-42 (2002).



\bibitem{Labbi-einstein} Labbi, M.-L., \emph{On compact manifolds with positive Einstein curvature},
Geometria Dedicata, 108, 205-217 (2004).

\bibitem{Labbi-double} Labbi, M.-L., \emph{ Double forms, curvature structures and the
$(p,q)$-curvatures}, Transactions of the American Mathematical Society, 357, n10, 3971-3992 (2005).

\bibitem{Lawson-Michelsohn} Lawson, H. B., Michelsohn,  M. L., \emph{Spin Geometry}, 
Princeton University Press (1989).

\bibitem{Maillot} Maillot, H., \emph{ Sur l'op\'erateur de courbure d'une vari\'et\'e riemannienne},
th\`ese, Universit\'e Claude Bernard, Lyon (1974).



\bibitem{Tachibana} Tachibana, S., \emph{The mean curvature for $p$-plane}, Journal of Differential Geometry, 8, 47-52 (1973).







\bibitem{Thorpe} Thorpe, J. A., \emph{Some remarks on the Gauss-Bonnet
integral}, Journal of Mathematics and Mechanics, Vol. 18, No. 8 (1969).

\end{thebibliography}
\end{document}